\documentclass[a4paper,12pt]{article}
\usepackage{enumerate}
\usepackage{latexsym}
\usepackage{amsfonts}
\usepackage[only,ninrm,elvrm,twlrm,sixrm,egtrm,tenrm]{rawfonts}
\usepackage{indentfirst}
\usepackage{amsmath}
\usepackage[noend]{algorithmic}
\usepackage{algorithm}
\usepackage{graphicx,psfrag}
\usepackage{graphics}
\usepackage{makeidx}
\newtheorem{thm}{Theorem}[section]
\newtheorem{cor}[thm]{Corollary}
\newtheorem{lem}[thm]{Lemma}

\newtheorem{exam}{Example}

\newtheorem{defn}{Definition}

\def\qed{\hfill \nopagebreak\rule{5pt}{8pt}}
\def\pf{\noindent {\it Proof.} }
\title{\bf NP-completeness of 4-incidence colorability of semi-cubic
graphs\footnote{Supported by NSFC.}}

\author{
\small Xueliang Li and Jianhua Tu\\
[3mm]
\small Center for Combinatorics and LPMC\\
\small Nankai University, Tianjin 300071, P.R. China\\
}
\date{}

\begin{document}

\maketitle
\begin{abstract}

The incidence coloring conjecture, proposed by Brualdi and Massey in
1993, states that the incidence coloring number of every graph is at
most ${\it \Delta}+2$, where ${\it \Delta}$ is the maximum degree of
a graph. The conjecture was shown to be false in general by Guiduli
in 1997, following the work of Algor and Alon. However, in 2005
Maydanskiy proved that the conjecture holds for any graph with ${\it
\Delta}\leq 3$. It is easily deduced that the incidence coloring
number of a semi-cubic graph is 4 or 5. In this paper, we show that
it is already NP-complete to determine if a semi-cubic graph is
4-incidence colorable, and therefore it is NP-complete to determine
if a general graph is $k$-incidence colorable.\\
[0.1in]{\bf Keywords:} incidence coloring number, $k$-incidence
colorable, strong-vertex coloring, semi-cubic graph, NP-complete.\\[0.2cm]

\end{abstract}

\section{Introduction}

In this paper we consider undirected, finite and simple graphs only,
and use standard notations in graph theory (see \cite{bon}). Let
$G=(V,E)$ be a graph, and let
$$I(G)=\{(v,e):v\in V,e\in E,\text{ and}\ v \text{ is incident with}\ e\}$$
be the set of all {\it incidences} of $G$. We say that two
incidences $(v,e)$ and $(w,f)$ are {\it adjacent} if one of the
following holds: (1) {\sl $v=w$,} (2) {\sl $e=f$ and} (3) {\sl the
edge $vw$ equals to $e$ or $f$.}

Following Shiu et al. \cite{shiu} we view $G$ as a digraph by
splitting each edge $uv$ into two opposite arcs $(u,v)$ and $(v,u)$.
For $e=uv$, we identify the incidence $(u,e)$ with the arc $(u,v)$.
By a slight abuse of notation we will refer to the incidence $(u,v)$
whenever it is convenient to do so. Two distinct incidences $(u,v)$
and $(x,y)$ are {\sl adjacent} if one of the following holds: (i)
{\sl $u=x$,} (ii) {\sl $u=y$ and $v=x$} (iii) {\sl $v=x$.} The
configurations associated with (i)-(iii) are pictured in Figure 1.

\begin{figure}
\begin{center}
\includegraphics[width=10cm]{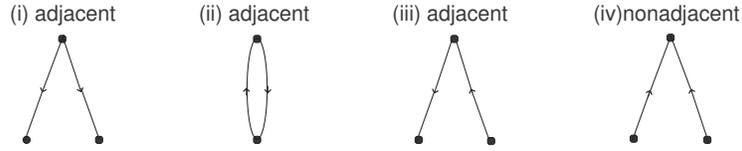}
\caption{Examples of adjacent and nonadjacent incidences.}
\end{center}
\end{figure}

We define an {\sl incidence coloring} of $G$ to be a coloring of its
incidences in which adjacent incidences are assigned different
colors. If $\sigma:I(G)\rightarrow S$ is an incidence coloring of
$G$ and $k=|S|$, then $G$ is called {\sl $k$-incidence colorable}
and $\sigma$ is a {\it $k$-incidence coloring}, where $S$ is a
color-set. The {\sl incidence coloring number} of $G$, denoted by
$\chi_{i}(G)$, is the smallest number of colors in an incidence
coloring.

This concept was first developed by Brualdi and Massey \cite{bru} in
1993. They proposed the {\it incidence coloring conjecture (ICC)},
which states that for every graph $G$, $\chi_{i}(G)\leq
{\it\Delta}+2$, where ${\it\Delta}={\it\Delta}(G)$ is the maximum
degree of $G$. In 1997 Guiduli \cite{gui} pointed out that the ICC
was solved in negative, following an example in \cite{Alon}. He
considered the Paley graphs of order $p$ with $p\equiv \ 1 \ (mod \
4)$. Following the analysis in \cite{Alon}, he showed that
$\chi_{i}(G)\geq{\it\Delta}+\Omega(log{\it\Delta})$, where
$\Omega=\frac{1}{8}-o(1)$. By using a tight upper bound for directed
star arboricity, he obtained the upper bound
$\chi_{i}(G)\leq{\it\Delta}+O(log{\it\Delta})$.

Brualdi and Massey \cite{bru} determined the incidence coloring
numbers of trees, complete graphs and complete bipartite graphs.
They also gave a simple bound for the incidence coloring number as
follows.

\begin{thm}(Brualdi and Massey \cite{bru})\label{bound}
For every graph $G$,
${\it\Delta}(G)+1\leq\chi_{i}(G)\leq2{\it\Delta}(G)$.
\end{thm}

In \cite{shiu} Shiu et al. proved that $\chi_{i}(G)\leq 5$ for
several classes of cubic (3-regular) 2-connected graphs $G$,
including Hamiltonian cubic graphs. In 2005 Maydanskiy \cite{may}
proved that the conjecture (ICC) holds for all graphs with
${\it\Delta}\leq 3$.

\begin{thm} (Maydanskiy \cite{may})\label{hold}
For every graph $G$ with ${\it\Delta}(G)\leq 3$, $\chi_{i}(G)\leq5$.
\end{thm}

\begin{defn}
For a graph $G$ with ${\it \Delta}=3$, if the degree of any vertex
of $G$ is 1 or 3, then the graph $G$ is called a semi-cubic graph.
\end{defn}

Using Theorems \ref{bound} and \ref{hold}, we have the following
corollary.

\begin{cor}
The incidence coloring number of a semi-cubic graph is 4 or 5.
\end{cor}

In this paper, we show the following result.

\begin{thm}\label{main}
It is NP-complete to determine if a semi-cubic graph is 4-incidence
colorable. Therefore, it is NP-complete to determine if a general
graph is $k$-incidence colorable, or in other words, it is
NP-complete to determine the incidence coloring number for a general
graph.
\end{thm}

\section{Incidence coloring of semi-cubic graphs}

In a given graph $G$, $N_{G}(v)$ denotes the set of vertices of $G$
adjacent to $v$, and $d_{G}(v)=|N_{G}(v)|$ is the degree of a vertex
$v$ in $G$. A vertex of degree $k$ is called a $k$-vertex. We denote
the set of all the incidences of the form $(v,u)$ and $(u,v)$ by
$O_{v}$ and $I_{v}$, respectively, where $u$ is adjacent with $v$.
Let $CO_{v}$ and $CI_{v}$ denote the set of colors assigned to
$O_{v}$ and $I_{v}$ in an incidence coloring of $G$, respectively.

We define a {\sl strong vertex coloring} of $G$ to be a proper
vertex coloring such that for any $u,w\in N_{G}(v)$, $u$ and $w$ are
assigned distinct colors. If $\sigma: V(G)\rightarrow S$ is a strong
vertex coloring of $G$ and $k=|S|$, then $G$ is called {\sl
$k$-strong-vertex colorable} and $\sigma$ is a {\it
$k$-strong-vertex coloring} of $G$, where $S$ is a color-set. And we
say that $G$ is $k$-strong-vertex colored.

\begin{lem}\label{strong}
Given a semi-cubic graph $G$, $G$ is 4-incidence colorable if and
only if $G$ is 4-strong-vertex colorable.
\end{lem}

\pf Since any two incidences in $O_{v}$ are adjacent, $|CO_{v}|$ is
equal to the degree of vertex $v$ in an incidence coloring of $G$.
Given a semi-cubic graph $G$, if $G$ is 4-incidence colorable and
$\sigma$ is a 4-incidence coloring, then $|CI_{v}|$ is 1 for every
vertex $v$ of $G$. We can color vertex $v$ using $CI_{v}$ and obtain
a vertex coloring. Since incidences $(u,v)$ and $(v,u)$ are
adjacent, $CI_{v}\neq CI_{u}$ and the vertex coloring is proper. If
$u,w\in N_{G}(v)$, incidences $(v,u)$ and $(v,w)$ are adjacent. Thus
$u$ and $w$ are assigned distinct colors in the vertex coloring. So
the vertex coloring is a 4-strong-vertex coloring and $G$ is
4-strong-vertex colorable.

Suppose there exists a 4-strong-vertex coloring $\sigma$ of $G$, we
color the incidence $(u,v)$ of $G$ by the color $c(v)$, where $c(v)$
denotes the color assigned to the vertex $v$ in $\sigma$. If two
incidences $(u,v)$ and $(x,y)$ are adjacent, then one of the
following holds: (1) $u=x$, $v\in N_{G}(u)$ and $y\in N_{G}(u)$; (2)
$u=y$ and $v=x$, $vy\in E(G)$; (3) $v=x$, $vy\in E(G)$. From the
definition of strong vertex coloring, incidences $(u,v)$ and $(x,y)$
are assigned different colors. So $G$ is 4-incidence colorable .\qed

\section{The blocks used in the construction}

For terminology and known results of NP-completeness we refer to the
book \cite{gra}. The $3SAT$ problem is stated as follows:\\

\noindent {\bf \underline {$3SAT$}}

\noindent INSTANCE: Set $U$ of variables, collection $\mathcal{C}$
of clauses over $U$ such that each clause $C_i\in \mathcal{C}$ has
$|C_i|=3$.

\noindent QUESTION: Is there a truth assignment for $U$ such that
every $C_i\in \mathcal{C}$ is true ?\\

It is clear that both the 4-incidence colorable problem for
semi-cubic graphs and the 4-strong-vertex colorable problem for
semi-cubic graphs are in the class NP. To prove their
NP-completeness, we exhibit a polynomial reduction from the known
NP-complete problem 3SAT. Given an instance $\mathcal{C}$ of the
problem 3SAT, we will show how to construct a semi-cubic graph $G$
of polynomial size in terms of the size of the instance
$\mathcal{C}$ such that $G$ is 4-strong-vertex colorable if and
only if $\mathcal{C}$ is satisfiable, which, from Lemma
\ref{strong}, implies that $G$ is 4-incidence colorable if and
only if $\mathcal{C}$ is satisfiable. The semi-cubic graph $G$
will be constructed from some pieces or ``blocks" which carry out
specific tasks. Information will be carried between blocks by
pairs of vertices. In a 4-strong-vertex coloring of the graph $G$,
such a pair of vertices is said to represent the value $T$
(``true") if the vertices have the same color, and to represent
$F$ (``false") if the vertices have distinct colors. In the
following we always use $S=\{1,2,3,4\}$ to denote the set of
colors.

\subsection{The switch blocks}

First, we give a special semi-cubic graph $H$ as shown in Figure 2.
We call $H$ a Kite graph. It is easy to check that if the Kite graph
$H$ is 4-strong-vertex colored, all 1-vertex $e, f, g$ have the same
color.

\begin{figure}
\begin{center}
\includegraphics[width=4cm]{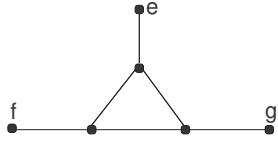}
\caption{The Kite graph $H$.}
\end{center}
\end{figure}

The switch block is shown with its symbol in Figure 3. It may be
checked that this block is 4-strong-vertex colorable. If this block
is 4-strong-vertex colored, one of the pairs of vertices marked
$a,b$ or $c,d$ must have the same color and the remaining pair of
vertices must have distinct colors. $c(l)$ is equal to $c(l')$ for
$l\in \{a,b,c,d\}$ in any 4-strong-vertex coloring of this block. In
fact, there are only two non-equivalence ways to color the two pairs
of vertices $a,b$ and $c,d$, which is shown in Figure 4.

\begin{figure}
\begin{center}
\includegraphics[width=8cm]{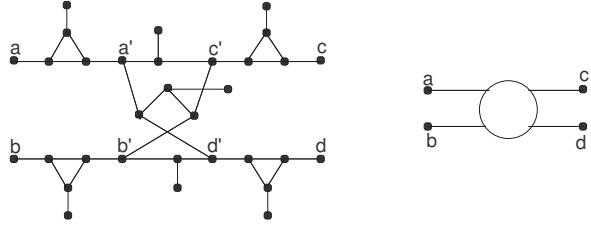}
\caption{The switch block and its symbolic representation.}
\end{center}
\end{figure}

\begin{figure}
\begin{center}
\includegraphics[width=12cm]{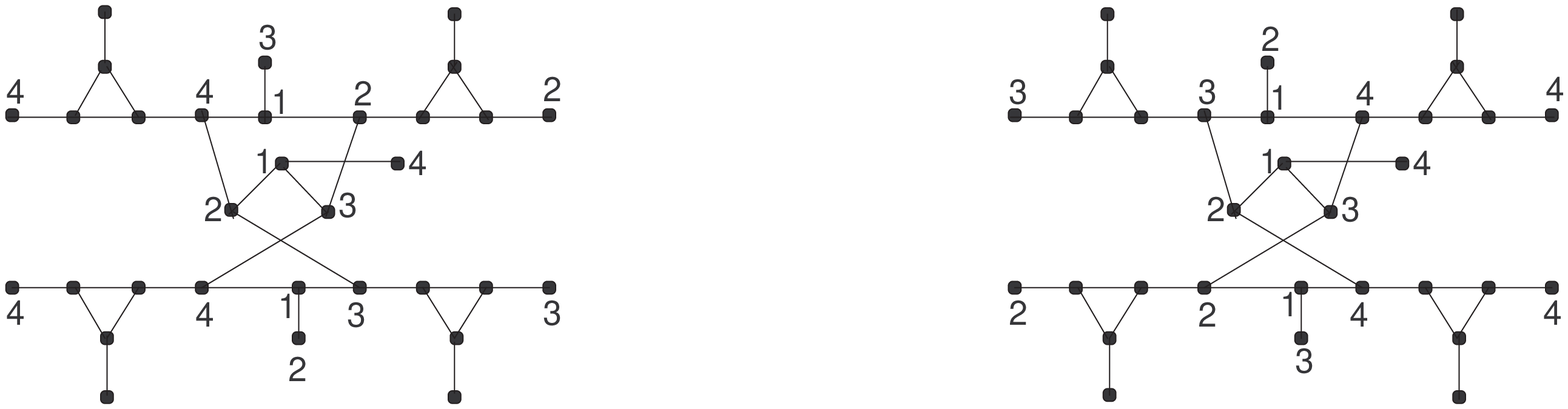}
\caption{Two non-equivalence ways to color the two pairs of
vertices $a,b$ and $c,d$.}
\end{center}
\end{figure}

Regarding the pair of vertices $a,b$ as the input and the pair
$c,d$ as the output, the block changes a representation of $T$ to
a representation of $F$, and vice versa.

\subsection{The variable blocks}

We construct the big switch block using three switch blocks by
identifying the output of the first switch block with the input of
the second switch block, identifying the output of the second switch
block with the input of the third switch block. It is shown with its
symbol in Figure 5. Regarding the pair of vertices $e,f$ as the
input and the pair $p,r$ as the output of the big switch block.

\begin{figure}
\begin{center}
\includegraphics[width=10cm]{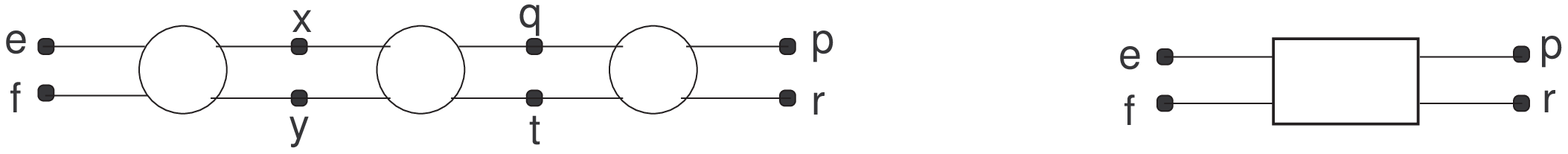}
\caption{The big switch block and its symbolic representation.}
\end{center}
\end{figure}

We state the following lemma.

\begin{lem}\label{big}
Given a big switch block $G$ for which the input and output are $e,f
$ and $p,r$, respectively, $G$ is 4-strong-vertex colorable. And if
$G$ is 4-strong-vertex colored, one of the following holds:
\begin{description}
\item[1] if $c(e)=c(f)$, then $c(p)\neq c(r)$ and $c(p)$,$c(r)$
may be any two different colors in $S=\{1,2,3,4\}$; \item[2] if
$c(e)\neq c(f)$, then $c(p)=c(r)$ and $c(p)$ may be any color in
$S=\{1,2,3,4\}$,
\end{description}
where $c(v)$ denotes the color assigned to the vertex $v$ of $G$.
\end{lem}

\pf Since the switch block is 4-strong-vertex colorable, it is easy
to check that the big switch block $G$ is also 4-strong-vertex
colorable.

If $G$ is 4-strong-vertex colored and $c(e)=c(f)$, without loss of
generality, we suppose that $c(e)$ is 1. The color-set $\{c(x),
c(y)\}$ of the output of the first switch block may be any two
different colors in $\{2,3,4\}$. For the output $t,q$ of the second
switch block, $c(q)$ is equal to $c(t)$ and may be any color in
$S=\{1, 2, 3, 4\}$. So, $c(p)\neq c(r)$ and $c(p)$, $c(r)$ may be
any two different colors in $S=\{1, 2, 3, 4\}$.

If $G$ is 4-strong-vertex colored and $c(e)\neq c(f)$, without loss
of generality, we suppose that $\{c(e),c(f)\}$ is $\{1,2\}$. Then
$c(x)$ is equal to $c(y)$ and may be any color in $\{3,4\}$. Thus
$c(q)\neq c(t)$, and $c(q)$ and $c(t)$ may be any two different
colors in $S=\{1,2,3,4\}$. So, $c(p)=c(r)$ and $c(p)$ may be any
color in $S=\{1,2,3,4\}$.\qed

From Lemma \ref{big}, the big switch block changes a representation
of $T$ to a representation of $F$, and vice versa.

The truth or falsity of each variable $u_{i}$ will be represented by
a variable block as shown in Figure 6, in which the blocks have,
respectively, $1,2,\cdots,6$ pairs of output vertices. In general,
the number of output pairs in the block representing $u_{i}$ should
be equal to the number of total appearances of $u_{i}$ or
$\overline{u_{i}}$ among the clauses of $\mathcal{C}$. If $k$ pairs
of output vertices is needed, we construct a variable block which is
made from $k$ big switch blocks and
$2\cdot\lfloor\frac{k-1}{2}\rfloor$ Kite graphs $H$.

\begin{figure}
\begin{center}
\includegraphics[width=8cm]{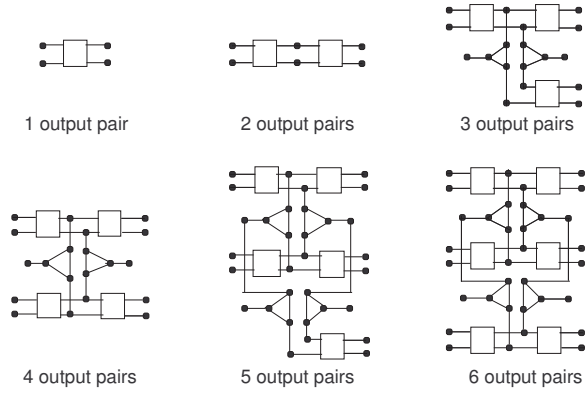}
\caption{The variable blocks having $1,2,\cdots,6$ output pairs of
vertices respectively. More generally it is made from $k$ big switch
blocks and $2\cdot\lfloor\frac{k-1}{2}\rfloor$ Kite graphs $H$ and
has $k$ output pairs.}
\end{center}
\end{figure}

From Lemma \ref{big} and the construction of the big switch block,
the following lemma is obvious.

\begin{lem}\label{var}
Given a variable block $G$, $G$ is 4-strong-vertex colorable. And in
any 4-strong-vertex coloring of $G$, all the output pairs must
represent the same value. If the output pairs represent $T$
(``true"), then the color-set of any output pairs may be any color
in $S=\{1,2,3,4\}$. If, on the other hand, the output pairs
represent $F$ (``false"), then the color-set of any output pairs may
be any two different colors in $S=\{1,2,3,4\}$. \qed
\end{lem}

\subsection{The clause blocks}

The truth of each clause $C_{j}$ will be tested by a clause block as
shown in Figure 7. The block is constructed from 3 switch blocks and
a cycle of length 10.

\begin{figure}
\begin{center}
\includegraphics[width=4cm]{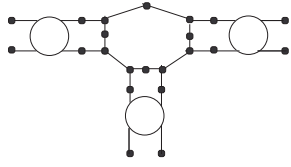}
\caption{The clause block made from 3 switch blocks and a cycle of
length 10 and having 3 input pairs of vertices.}
\end{center}
\end{figure}

The following lemma is crucial for proving our main theorem.

\begin{lem}\label{sat}
The clause block is 4-strong-vertex colorable if and only if the
three input pairs of vertices do not all represent $F$.
\end{lem}

\pf Given a clause block $G$ as shown in Figure 8, we suppose that
the 3 output pairs of vertices are $\{r,r'\}$, $\{w,w'\}$ and
$\{z,z'\}$, respectively.

We suppose that the three input pairs of vertices all represent $F$.
If $G$ is 4-strong-vertex colorable and $\sigma$ is a
4-strong-vertex coloring, then $c(r)=c(r')$, $c(w)=c(w')$ and
$c(z)=c(z')$, where $c(v)$ denotes the color assigned to the vertex
$v$ in $\sigma$. Without loss of generality, we suppose that $c(z)$,
$c(b)$ and $c(d)$ are 1, 2 and 3, respectively. Then $c(a)=c(e)=4$,
$c(f)=2$. Since $c(w)$ is equal to $c(w')$, $\{c(w),c(g)\}=\{c(w'),
c(g)\}=\{1,3\}$. Thus $c(h)=4$, $c(i)=2$. By the same token, $c(p)$
must be 4. So $c(a)=c(p)=4$ and $pa\in E(G)$, a contradiction.

\begin{figure}
\begin{center}
\includegraphics[width=10cm]{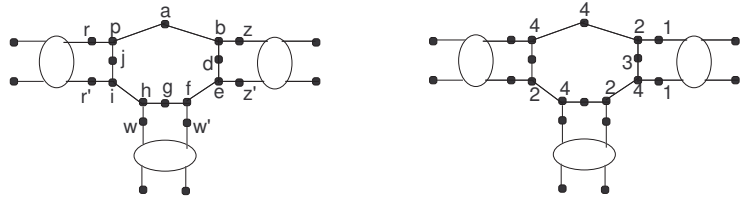}
\caption{The three input pairs of vertices all represent $F$.}
\end{center}
\end{figure}

If the three input pairs of vertices do not all represent $F$, we
can give a 4-strong-vertex coloring of $G$, which is shown in Figure
9. Thus $G$ is 4-strong-vertex colorable.

\begin{figure}
\begin{center}
\includegraphics[width=9cm]{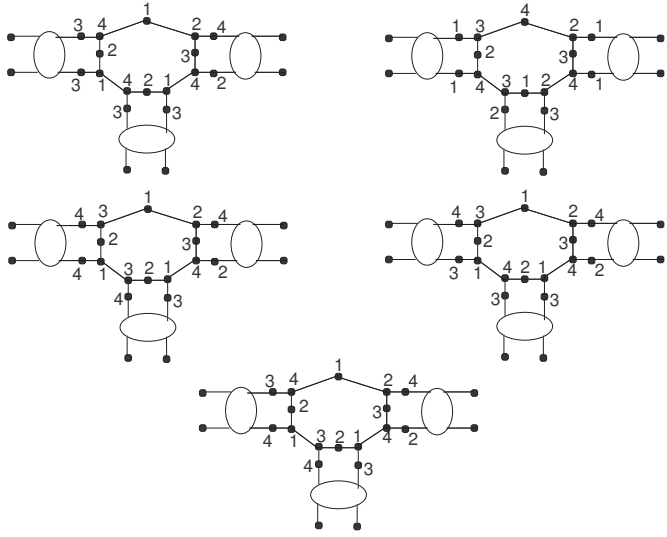}
\caption{The three input pairs of vertices do not all represent
$F$.}
\end{center}
\end{figure}

So the clause block is 4-strong-vertex colorable if and only if the
three input pairs of vertices do not all represent $F$.\qed

\section{Main result}

In this section, we prove the main result that it is NP-complete to
determine whether the incidence coloring number of a semi-cubic
graph is 4 or 5. Thus this problem has no polynomial time algorithm
unless $P=NP$.

{\bf{\it Proof of Theorem \ref{main}.}} The problem is clearly in
the class NP. We exhibit a polynomial reduction from the problem
3SAT. Consider an instance $\mathcal{C}$ of 3SAT and construct from
it a semi-cubic graph $G$ as follows.

Each variable $u_{i}$ corresponds to a variable block $U_{i}$ with
one output pair of vertices associated with each appearance of
$u_{i}$ or $\overline{u_{i}}$ among the clauses of $\mathcal{C}$.
Each clause $C_{j}$ corresponds to a clause block $B_{j}$. Suppose
literal $x_{j,k}$ in clause $C_{j}$ is the variable $u_{i}$. Then
identify the $k$-th input pair of $B_{j}$ with the associate output
pair of $U_{i}$. If, on the other hand, $x_{j,k}$ is
$\overline{u_{i}}$, then insert an switch block between the $k$-th
input pair of $B_{j}$ and the associated output pair of $U_{i}$. The
resulting graph $G'$ has some 2-vertices. For every 2-vertex $u$, we
add a pendant edge $(u,u')$. Now the resulting graph $G$ is a
semi-cubic graph.

We shall show that the incidence coloring number of the semi-cubic
graph $G$ is 4 if and only if there is a truth assignment to the
variables which simultaneously satisfies all the clauses in
$\mathcal{C}$.

First, we suppose that the incidence coloring number of $G$ is 4.
From Lemma \ref{strong}, $G$ is 4-strong-vertex colorable. Given a
4-strong-vertex coloring of $G$, we assign the value to each
variable $u_{i}$ as the one that all the output pairs of the
variable block corresponding to $u_{i}$ represent. By Lemma
\ref{sat}, the assignment to the variables simultaneously satisfies
all the clauses in $\mathcal{C}$.

Now, we suppose that there is a truth assignment to the variables
which simultaneously satisfies all the clauses in $\mathcal{C}$. For
each variable $u_{i}$, we can 4-strong-vertex color the variable
block corresponding to $u_{i}$ such that the value represented by
the output pairs is equal to the value in the truth assignment. And
by Lemma \ref{var}, the color-set of any output pairs may be any
color or any two distinct colors in $S=\{1,2,3,4\}$. So by Lemma
\ref{sat}, $G$ is 4-strong-vertex colorable. Finally, by Lemma
\ref{strong}, the incidence coloring number of $G$ is 4.

It is clear that the semi-cubic graph $G$ is constructed from $3m$
big switch blocks, $m$ clause blocks, at most $2\cdot \lfloor
\frac{3m-1}{2}\rfloor$ Kite graphs, and at most $3m$ switch blocks
and some pendant edges. Thus, $|V(G)|$ is at most $521m$, and
$|E(G)|$ is at most $626m$, and hence the semi-cubic graph $G$ can
be constructed from $\mathcal{C}$ polynomially in $m$, the size of
$\mathcal{C}$. The proof is now complete.\qed

\begin{exam}
Let $\mathcal{C}=\{C_{1},C_{2}\}$ an instance of the problem 3SAT
and
$$C_{1}=u_{1}\vee\overline{u_{2}}\vee u_{3},$$ $$C_{2}=u_{2}\vee
u_{3}\vee \overline{u_{4}}.$$ By the proof of Theorem \ref{main}, we
can construct the graph $G'$ as shown in Figure 10 and obtain the
semi-cubic $G$ by adding a pendant edge to every 2-vertex of $G'$.

\begin{figure}
\begin{center}
\includegraphics[width=7cm]{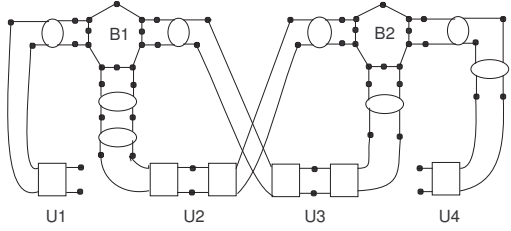}
\caption{The graph $G'$. We can obtain a semi-cubic graph $G$ by
adding a pendant edge to every 2-vertex of $G'$.}
\end{center}
\end{figure}
\end{exam}

\end{document}